\pdfoutput=1

\documentclass[11pt,letterpaper]{article}
\usepackage{naaclhlt2013}
\usepackage{times}
\usepackage{latexsym}
\usepackage{qtree}

\usepackage{amssymb}
\usepackage{amsmath}

\usepackage{epsfig}

\usepackage{color}

\usepackage{txfonts}

\newcommand{\semantics}[1]{[\![ #1 ]\!]}

\title{Towards a Formal Distributional Semantics:\\Simulating Logical Calculi with Tensors}

\author{Edward Grefenstette\\
	    University of Oxford\\ Department of Computer Science\\
	    Wolfson Building, Parks Road\\
	    Oxford OX1 3QD, UK\\
	    {\tt edward.grefenstette@cs.ox.ac.uk}
	  }

\date{}

\begin{document}
\maketitle
\begin{abstract}
  The development of compositional distributional models of semantics reconciling the empirical aspects of distributional semantics with the compositional aspects of formal semantics is a popular topic in the contemporary literature. This paper seeks to bring this reconciliation one step further by showing how the mathematical constructs commonly used in compositional distributional models, such as tensors and matrices, can be used to simulate different aspects of predicate logic. 

  This paper discusses how the canonical isomorphism between tensors and multilinear maps can be exploited to simulate a full-blown quantifier-free predicate calculus using tensors. It provides tensor interpretations of the set of logical connectives required to model propositional calculi. It suggests a variant of these tensor calculi capable of modelling quantifiers, using few non-linear operations. It finally discusses the relation between these variants, and how this relation should constitute the subject of future work.
\end{abstract}

\section{Introduction} 
\label{sec:introduction}

The topic of compositional distributional semantics has been growing in popularity over the past few years. This emerging sub-field of natural language semantic modelling seeks to combine two seemingly orthogonal approaches to modelling the meaning of words and sentences, namely \emph{formal semantics} and \emph{distributional semantics}.

These approaches, summarised in Section~\ref{sec:background}, differ in that formal semantics, on the one hand, provides a neatly compositional picture of natural language meaning, reducing sentences to logical representations; one the other hand, distributional semantics accounts for the ever-present ambiguity and polysemy of words of natural language, and provides tractable ways of learning and comparing word meanings based on corpus data.

Recent efforts, some of which are briefly reported below, have been made to unify both of these approaches to language modelling to produce \emph{compositional distributional models of semantics}, leveraging the learning mechanisms of distributional semantics, and providing syntax-sensitive operations for the production of representations of sentence meaning obtained through combination of corpus-inferred word meanings. These efforts have been met with some success in evaluations such as phrase similarity tasks~\cite{mitchell2008vector,mitchell2009language,Grefenstette2011a,Kartsaklis2012}, sentiment prediction~\cite{socherEMNLP12}, and paraphrase detection~\cite{blacoe70comparison}.

While these developments are promising with regard to the goal of obtaining learnable-yet-structured sentence-level representations of language meaning, part of the motivation for unifying formal and distributional models of semantics has been lost. The compositional aspects of formal semantics are combined with the corpus-based empirical aspects of distributional semantics in such models, yet the logical aspects are not. But it is these logical aspects which are so appealing in formal semantic models, and therefore it would be desirable to replicate the inferential powers of logic within compositional distributional models of semantics.

In this paper, I make steps towards addressing this lost connection with logic in compositional distributional semantics. In Section~\ref{sec:background}, I provide a brief overview of formal and distributional semantic models of meaning. In Section~\ref{sec:tensors_as_multilinear_maps}, I give mathematical foundations for the rest of the paper by introducing tensors and tensor contraction as a way of modelling multilinear functions. In Section~\ref{sec:tensor_based_model_theoretic_predicate_calculi}, I discuss how predicates, relations, and logical atoms of a quantifier-free predicate calculus can be modelled with tensors. In Section~\ref{sec:logical_connectives_with_tensors}, I present tensorial representations of logical operations for a complete propositional calculus. In Section~\ref{sec:quantifiers_and_non_linearity}, I discuss a variant of the predicate calculus from Section~\ref{sec:tensor_based_model_theoretic_predicate_calculi} aimed at modelling quantifiers within such tensor-based logics, and the limits of compositional formalisms based only on multilinear maps. I conclude, in Section~\ref{sec:conclusions_and_future_work}, by suggesting directions for further work based on the contents of this paper.

This paper does not seek to address the question of how to determine how words should be translated into predicates and relations in the first place, but rather shows how such predicates and relations can be modelled using multilinear algebra. As such, it can be seen as a general theoretical contribution which is independent from the approaches to compositional distributional semantics it can be applied to. It is directly compatible with the efforts of~\newcite{Coecke2010} and~\newcite{grefenstette2013multistep}, discussed below, but is also relevant to any other approach making use of tensors or matrices to encode semantic relations.


\section{Related work} 
\label{sec:background}

\begin{figure}[ht!]
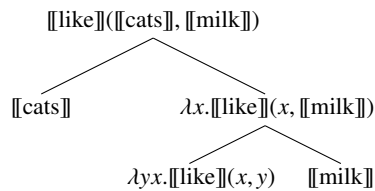

\begin{center}
    {\footnotesize       
            \begin{tabular}{l|l}
                Syntactic Analysis & Semantic Interpretation\\
                \hline\\[-0.3cm]
                S $\Rightarrow$ NP VP & $\semantics{VP}(\semantics{NP})$\\
                NP $\Rightarrow$ cats, milk, etc. & $\semantics{{\textrm{cats}}},\,\semantics{\textrm{milk}},\,\ldots$\\
                VP $\Rightarrow$ Vt NP & $\semantics{Vt}(\semantics{NP})$\\
                Vt $\Rightarrow$ like, hug, etc. & $\lambda yx.\semantics{\textrm{like}}(x,y),\,\ldots$ \\
            \end{tabular}\\[.5cm]

            \Tree [.$\semantics{\textrm{like}}(\semantics{\textrm{cats}},\semantics{\textrm{milk}})$ $\semantics{\textrm{cats}}$ [.$\quad \lambda x.\semantics{\textrm{like}}(x,\semantics{\textrm{milk}})$  $\quad\lambda yx.\semantics{\textrm{like}}(x,y)$  $\semantics{\textrm{milk}}$ ] ]
            }
\end{center}
\caption{A simple formal semantic model.}
\label{fig:simpleformalsemantics}
\end{figure}

Formal semantics, from the Montagovian school of thought \cite{Montague1974,Dowty1981montague}, treats natural languages as programming languages which compile down to some formal language such as a predicate calculus. The syntax of natural languages, in the form of a grammar, is augmented by semantic interpretations, in the form of expressions from a higher order logic such as the lambda-beta calculus. The parse of a sentence then determines the combinations of lambda-expressions, the reduction of which yields a well-formed formula of a predicate calculus, corresponding to the semantic representation of the sentence. A simple formal semantic model is illustrated in Figure~\ref{fig:simpleformalsemantics}.

Formal semantic models are incredibly powerful, in that the resulting logical representations of sentences can be fed to automated theorem provers to perform textual inference, consistency verification, question answering, and a host of other tasks which are well developed in the literature (e.g.~see~\cite{loveland1978automated} and~\cite{fitting1996first}). However, the sophistication of such formal semantic models comes at a cost: the complex set of rules allowing for the logical interpretation of text must either be provided \emph{a priori}, or learned. Learning such representations is a complex task, the difficulty of which is compounded by issues of ambiguity and polysemy which are pervasive in natural languages.

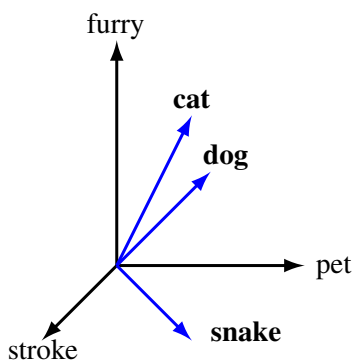
\begin{figure}
  \begin{center}
    \setlength{\unitlength}{0.5mm}
            \begin{picture}(100, 100)
              \thicklines
              \put(22,82){\text{furry}}
              \linethickness{0.4mm}
              \put(30, 20){\vector(1, 0){50}}
              \put(1,-5) {\text{stroke}}
              \put(30, 20){\vector(0, 1){60}}
              \put(83,18){\text{pet}}
              \put(30, 20){\vector(-1, -1){20}}
                \thinlines
                  \linethickness{0.4mm}
              \put(30, 20){\begin{color}{blue}\vector(1, 2){20}\end{color}}
              \put(30, 20){\begin{color}{blue}\vector(1, 1){25}\end{color}}
              \put(30, 20){\begin{color}{blue}\vector(0.05, -0.05){20}\end{color}}
              \put(45,62){\text{{\bf cat}}}
              \put(53,47){\text{\bf dog}}
              \put(55,0){\text{\bf snake}}
            \end{picture}
  \end{center}
  \caption{A simple distributional semantic model.}
  \label{fig:simpledistsemantics}
\end{figure}

In contrast, distributional semantic models, best summarised by the dictum of \newcite{Firth1957} that ``You shall know a word by the company it keeps,'' provide an elegant and tractable way of learning semantic representations of words from text. Word meanings are modelled as high-dimensional vectors in large semantic vector spaces, the basis elements of which correspond to contextual features such as other words from a lexicon. Semantic vectors for words are built by counting how many time a target word occurs within a context (e.g.~within $k$ words of select words from the lexicon). These context counts are then normalised by a term frequency-inverse document frequency-like measure (e.g.~TF-IDF, pointwise mutual information, ratio of probabilities), and are set as the basis weights of the vector representation of the word's meaning. Word vectors can then be compared using geometric distance metrics such as cosine similarity, allowing us to determine the similarity of words, cluster semantically related words, and so on. Excellent overviews of distributional semantic models are provided by~\newcite{Curran2004} and~\newcite{mitchell2011composition}. A simple distributional semantic model showing the spacial representation of words `dog', `cat' and `snake' within the context of feature words `pet', `furry', and `stroke' is shown in Figure~\ref{fig:simpledistsemantics}.

Distributional semantic models have been successfully applied to tasks such as word-sense discrimination~\cite{schutze1998automatic}, thesaurus extraction~\cite{Grefenstette1994}, and automated essay marking~\cite{Landauer1997}. However, while such models provide tractable ways of learning and comparing word meanings, they do not naturally scale beyond word length. As recently pointed out by~\newcite{turney2012domain}, treating larger segments of texts as lexical units and learning their representations distributionally (the `holistic approach') violates the principle of linguistic creativity, according to which we can formulate and understand phrases which we've never observed before, provided we know the meaning of their parts and how they are combined. As such, distributional semantics makes no effort to account for the compositional nature of language like formal semantics does, and ignores issues relating to syntactic and relational aspects of language.

Several proposals have been put forth over the last few years to provide vector composition functions for distributional models in order to introduce compositionality, thereby replicating some of the aspects of formal semantics while preserving learnability. Simple operations such as vector addition and multiplication, with or without scalar or matrix weights (to take word order or basic relational aspects into account), have been suggested \cite{zanzotto2010estimating,mitchell2008vector,mitchell2009language}. 

\newcite{smolensky1990tensor} suggests using the tensor product of word vectors to produce representations that grow with sentence complexity. \newcite{Clark2006} extend this approach by including basis vectors standing for dependency relations into tensor product-based representations. Both of these tensor product-based approaches run into dimensionality problems as representations of sentence meaning for sentences of different lengths or grammatical structure do not live in the same space, and thus cannot directly be compared. \newcite{Coecke2010} develop a framework using category theory, solving this dimensionality problem of tensor-based models by projecting tensored vectors for sentences into a unique vector space for sentences, using functions dynamically generated by the syntactic structure of the sentences. In presenting their framework, which partly inspired this paper, they describe how a verb can be treated as a logical relation using tensors in order to evaluate the truth value of a simple sentence, as well as how negation can be modelled using matrices.

A related approach, by \newcite{Baroni2010}, represents unary relations such as adjectives as matrices learned by linear regression from corpus data, and models adjective-noun composition as matrix-vector multiplication. \newcite{grefenstette2013multistep} generalise this approach to relations of any arity and relate it to the framework of \newcite{Coecke2010} using a tensor-based approach to formal semantic modelling similar to that presented in this paper.

Finally, \newcite{socherEMNLP12} apply deep learning techniques to model syntax-sensitive vector composition using non-linear operations, effectively turning parse trees into multi-stage neural networks. Socher shows that the non-linear activation function used in such a neural network can be tailored to replicate the behaviour of basic logical connectives such as conjunction and negation.


\section{Tensors and multilinear maps} 
\label{sec:tensors_as_multilinear_maps}

Tensors are the mathematical objects dealt with in multilinear algebra just as vectors and matrices are the objects dealt with in linear algebra. In fact, tensors can be seen as generalisations of vectors and matrices by introducing the notion of \emph{tensor rank}. Let the rank of a tensor be the number of indices required to describe a vector/matrix-like object in sum notation. A vector $\mathbf{v}$ in a space $V$ with basis $\{\mathbf{b}^V_i\}_i$ can be written as the weighted sum of the basis vectors:
\[
\mathbf{v} = \sum_{i}{c^v_i \mathbf{b}^V_i}
\]
where the $c^v_i$ elements are the scalar basis weights of the vector. Being fully described with one index, vectors are rank 1 tensors. Similarly, a matrix $\mathbf{M}$ is an element of a space $V \otimes W$ with basis $\{(\mathbf{b}^V_i,\mathbf{b}^W_j)\}_{ij}$ (such pairs of basis vectors of $V$ and $W$ are commonly written as $\{\mathbf{b}^V_i \otimes \mathbf{b}^W_j\}_{ij}$ in multilinear algebra). Such matrices are rank 2 tensors, as they can be fully described using two indices (one for rows, one for columns):
\[
\mathbf{M} = \sum_{ij}{c^M_{ij} \mathbf{b}^V_i \otimes \mathbf{b}^W_j}
\]
where the scalar weights $c^M_{ij}$ are just the $ij$th elements of the matrix.

A tensor $\mathbf{T}$ of rank $k$ is just a geometric object with a higher rank. Let $T$ be a member of $V_1 \otimes \ldots \otimes V_k$; we can express $T$ as follows, using $k$ indices $\alpha_1 \ldots \alpha_k$:
\[
\mathbf{T} = \sum_{\alpha_1 \ldots \alpha_k}{c^T_{\alpha_1 \ldots \alpha_k} \mathbf{b}^{V_1}_{\alpha_1} \otimes \ldots \otimes \mathbf{b}^{V_k}_{\alpha_k}}
\]
In this paper, we will be dealing with tensors of rank 1 (vectors), rank 2 (matrices) and rank 3, which can be pictured as cuboids (or a matrix of matrices).

Tensor contraction is an operation which allows us to take two tensors and produce a third. It is a generalisation of inner products and matrix multiplication to tensors of higher ranks. Let $\mathbf{T}$ be a tensor in $V_1 \otimes \ldots \otimes V_j \otimes V_k$ and $\mathbf{U}$ be a tensor in $V_k \otimes V_m \otimes \ldots \otimes V_n$. The contraction of these tensors, written $\mathbf{T} \times \mathbf{U}$, corresponds to the following calculation:
{\small
\begin{align*}
& \mathbf{T} \times \mathbf{U} = \\
& \quad \sum_{\alpha_1 \ldots \alpha_n}{c^T_{\alpha_1 \ldots\alpha_k}c^U_{\alpha_k \ldots \alpha_n} \mathbf{b}^{V_1}_{\alpha_1} \otimes \ldots \otimes \mathbf{b}^{V_j}_{\alpha_j} \otimes \mathbf{b}^{V_m}_{\alpha_m} \otimes \ldots \otimes \mathbf{b}^{V_n}_{\alpha_n}}
\end{align*}}
Tensor contraction takes a tensor of rank $k$ and a tensor of rank $n-k+1$ and produces a tensor of rank $n-1$, corresponding to the sum of the ranks of the input tensors minus $2$. The tensors must satisfy the following restriction: the left tensor must have a rightmost index spanning the same number of dimensions as the leftmost index of the right tensor. This is similar to the restriction that a $m$ by $n$ matrix can only be multiplied with a $p$ by $q$ matrix if $n = p$, i.e.~if the index spanning the columns of the first matrix covers the same number of columns as the index spanning the rows of the second matrix covers rows. Similarly to how the columns of one matrix `merge' with the rows of another to produce a third matrix, the part of the first tensor spanned by the index $k$ merges with the part of the second tensor spanned by $k$ by `summing through' the shared basis elements $\mathbf{b}^{V_k}_{\alpha_k}$ of each tensor. Each tensor therefore loses a rank while being joined, explaining how the tensor produced by $\mathbf{T} \times \mathbf{U}$ is of rank $k+(n-k+1) -2 = n-1$.

There exists an isomorphism between tensors and multilinear maps \cite{Bourbaki:1989,Lee:1997}, such that any curried multilinear map
\[
f : V_1 \to \ldots \to V_j \to V_k
\]
can be represented as a tensor $\mathbf{T}^f \in V_k \otimes V_j \otimes \ldots \otimes V_1$ (note the reversed order of the vector spaces), with tensor contraction acting as function application. This isomorphism guarantees that there exists such a tensor $\mathbf{T}^f$ for every $f$, such that the following equality holds for any $\mathbf{v}_1 \in V_1,\, \ldots,\, \mathbf{v}_j \in V_j$:
\[
f \mathbf{v}_1 \ldots \mathbf{v}_j = \mathbf{v}_k = \mathbf{T}^f \times \mathbf{v}_1 \times \ldots \times \mathbf{v}_j 
\]


\section{Tensor-based predicate calculi} 
\label{sec:tensor_based_model_theoretic_predicate_calculi}

In this section, I discuss how the isomorphism between multilinear maps and tensors described above can be used to model predicates, relations, and logical atoms of a predicate calculus. The four aspects of a predicate calculus we must replicate here using tensors are as follows: truth values, the logical domain and its elements (logical atoms), predicates, and relations. I will discuss logical connectives in the next section.

Both truth values and domain objects are the basic elements of a predicate calculus, and therefore it makes sense to model them as vectors rather than higher rank tensors, which I will reserve for relations. We first must consider the vector space used to model the boolean truth values of $\mathbb{B}$. \newcite{Coecke2010} suggest, as boolean vector space, the space $B$ with the basis $\{\top,\bot\}$, where $\top = [1\ 0]^{\top}$ is interpreted as `true', and $\bot = [0\ 1]^{\top}$ as `false'.

I assign to the domain $\mathcal{D}$, the set of objects in our logic, a vector space $D$ on $\mathbb{R}^{|\mathcal{D}|}$ with basis vectors $\{\mathbf{d}_i\}_i$ which are in bijective correspondence with elements of $\mathcal{D}$. An element of $\mathcal{D}$ is therefore represented as a one-hot vector in $D$, the single non-null value of which is the weight for the basis vector mapped to that element of $\mathcal{D}$. Similarly, a subset of $\mathcal{D}$ is a vector of $D$ where those elements of $\mathcal{D}$ in the subset have $1$ as their corresponding basis weights in the vector, and those not in the subset have $0$. Therefore there is a one-to-one correspondence between the vectors in $D$ and the elements of the power set $\mathcal{P}(\mathcal{D})$, provided the basis weights of the vectors are restricted to one of $0$ or $1$.

Each unary predicate $P$ in the logic is represented in the logical model as a set $M_P \subseteq \mathcal{D}$ containing the elements of the domain for which the predicate is true. Predicates can be viewed as a unary function $f_P : \mathcal{D} \to \mathbb{B}$ where 
\[
    f_P(x) =
    \left\{
    \begin{tabular}{l}
     $\top\quad$ if $x \in M_P$\\
     $\bot\quad$ otherwise
    \end{tabular}
    \right.
  \]
 These predicate functions can be modelled as rank 2 tensors in $B \otimes D$, i.e.~matrices. Such a matrix $\mathbf{M}^P$ is expressed in sum notation as follows:
\[
\mathbf{M}^P = \left( \sum_i{c^{M^P}_{1i} \top \otimes \mathbf{d}_i} \right) + \left(\sum_i{c^{M^P}_{2i} \bot \otimes \mathbf{d}_i} \right)
\]
The basis weights are defined in terms of the set $M_P$ as follows: $c^{M^P}_{1i} = 1$ if the logical atom $x_i$ associated with basis weight $\mathbf{d}_i$ is in $M_P$, and $0$ otherwise; conversely, $c^{M^P}_{2i} = 1$ if the logical atom $x_i$ associated with basis weight $\mathbf{d}_i$ is \emph{not} in $M_P$, and $0$ otherwise.

To give a simple example, let's consider a domain with three individuals, represented as the following one-hot vectors in $D$: $\mathbf{john} = [1\ 0\ 0]^\top$, $\mathbf{chris} = [0\ 1\ 0]^\top$, and $\mathbf{tom} = [0\ 0\ 1]^\top$. Let's imagine that Chris and John are mathematicians, but Tom is not. The predicate $P$ for `is a mathematician' therefore is represented model-theoretically as the set $M_P = \{chris,\, john\}$. Translating this into a matrix gives the following tensor for $P$:
\[
\mathbf{M}^P =
\left[
\begin{tabular}{ccc}
1 &  1 & 0\\
0 & 0 & 1 
\end{tabular}
\right]
\]
To compute the truth value of `John is a mathematician', we perform predicate-argument application as tensor contraction (matrix-vector multiplication, in this case):
\[
\mathbf{M}^P \times \mathbf{john} =
\left[
\begin{tabular}{ccc}
1 &  1 & 0\\
0 & 0 & 1 
\end{tabular}
\right]
\left[
\begin{tabular}{c}
0 \\
1 \\
0 \\
\end{tabular}
\right]
=
\left[
\begin{tabular}{c}
1 \\
0 \\
\end{tabular}
\right] = \top
\]
Likewise for `Tom is a mathematician':
\[
\mathbf{M}^P \times \mathbf{tom} =
\left[
\begin{tabular}{ccc}
1 &  1 & 0\\
0 & 0 & 1 
\end{tabular}
\right]
\left[
\begin{tabular}{c}
0 \\
0 \\
1 \\
\end{tabular}
\right]
=
\left[
\begin{tabular}{c}
0 \\
1 \\
\end{tabular}
\right] = \bot
\]

Model theory for predicate calculus represents any $n$-ary relation $R$, such as a verb, as the set $M_R$ of $n$-tuples of elements from $\mathcal{D}$ for which $R$ holds. Therefore such relations can be viewed as functions $f_R: \mathcal{D}^n \to \mathbb{B}$ where:
    \[
    f_R(x_1,\ldots,x_n) =
    \left\{
    \begin{tabular}{l}
    $\top\quad$ if $(x_1,\ldots,x_n) \in M_R$\\
    $\bot\quad$ otherwise
    \end{tabular}
    \right.
    \]
We can represent the boolean function for such a relation $R$ as a tensor $\mathbf{T}^R$ in $B \otimes \underbrace{D \otimes \ldots \otimes D}_n$:
\begin{align*}
& \mathbf{T}^R = & \left(\sum_{\alpha_1 \ldots \alpha_n}{c^{T^R}_{1\alpha_1 \ldots \alpha_n} \top \otimes \mathbf{d}_{\alpha_1} \otimes \ldots \otimes \mathbf{d}_{\alpha_n}}\right)\\
& & + \left(\sum_{\alpha_1 \ldots \alpha_n}{c^{T^R}_{2\alpha_1 \ldots \alpha_n} \bot \otimes \mathbf{d}_{\alpha_1} \otimes \ldots \otimes \mathbf{d}_{\alpha_n}}\right)
\end{align*}
As was the case for predicates, the weights for relational tensors are defined in terms of the set modelling the relation: $c^{T^R}_{1\alpha_1 \ldots \alpha_n}$ is $1$ if the tuple $(x,\, \ldots,\, z)$ associated with the basis vectors $\mathbf{d}_{\alpha_n} \ldots \mathbf{d}_{\alpha_1}$ (again, note the reverse order) is in $M_R$ and $0$ otherwise; and $c^{T^R}_{2\alpha_1 \ldots \alpha_n}$ is $1$ if the tuple $(x,\, \ldots,\, z)$ associated with the basis vectors $\mathbf{d}_{\alpha_n} \ldots \mathbf{d}_{\alpha_1}$ is \emph{not} in $M_R$ and $0$ otherwise.

To give an example involving relations, let our domain be the individuals John ($j$) and Mary ($m$). Mary loves John and herself, but John only loves himself. The logical model for this scenario is as follows:
\[
\mathcal{D} = \{j,m\} \qquad  M_\text{loves} = \{ (j,j),\, (m,m),\, (m,j) \}
\]
Distributionally speaking, the elements of the domain will be mapped to the following one-hot vectors in some two-dimensional space $D$ as follows: $\mathbf{j} = [1\ 0]^\top$ and $\mathbf{m} = [0\ 1]^\top$. The tensor for `loves' can be written as follows, ignoring basis elements with null-valued basis weights, and using the distributivity of the tensor product over addition:
\begin{align*}
& \mathbf{T}^{\text{loves}} =  \top \otimes ((\mathbf{d}_1 \otimes \mathbf{d}_1) + (\mathbf{d}_2 \otimes \mathbf{d}_2) + (\mathbf{d}_1 \otimes \mathbf{d}_2))\\
& \quad \quad\ \,  + (\bot \otimes \mathbf{d}_2 \otimes \mathbf{d}_1)
\end{align*}
Computing ``Mary loves John'' would correspond to the following calculation:
\begin{align*}
  & (\mathbf{T}^{\textrm{loves}} \times \mathbf{m}) \times \mathbf{j} =\\
  & \quad ((\top \otimes \mathbf{d}_2) + (\top \otimes \mathbf{d}_1)) \times \mathbf{j} = \top
\end{align*}
whereas ``John loves Mary'' would correspond to the following calculation:
\begin{align*}
  & (\mathbf{T}^{\textrm{loves}} \times \mathbf{j}) \times \mathbf{m} =\\
  & \quad ((\top \otimes \mathbf{d}_1) + (\bot \otimes \mathbf{d}_2)) \times \mathbf{m} = \bot
\end{align*}


\section{Logical connectives with tensors} 
\label{sec:logical_connectives_with_tensors}

In this section, I discuss how the boolean connectives of a propositional calculus can be modelled using tensors. Combined with the predicate and relation representations discussed above, these form a complete quantifier-free predicate calculus based on tensors and tensor contraction.

Negation has already been shown to be modelled in the boolean space described earlier by~\newcite{Coecke2010} as the swap matrix:
\[
\mathbf{T}^{\lnot} =
\left[
\begin{tabular}{cc}
0 & 1\\
1 & 0
\end{tabular}
\right]
\]
This can easily be verified:
\begin{align*}
  & \mathbf{T}^{\lnot} \times {\top} =
\left[
\begin{tabular}{cc}
0 & 1\\
1 & 0
\end{tabular}
\right]
\left[
\begin{tabular}{c}
1\\
0
\end{tabular}
\right]
= 
\left[
\begin{tabular}{c}
0\\
1
\end{tabular}
\right]
= {\bot}\\
& \mathbf{T}^{\lnot} \times {\bot} =
\left[
\begin{tabular}{cc}
0 & 1\\
1 & 0
\end{tabular}
\right]
\left[
\begin{tabular}{c}
0\\
1
\end{tabular}
\right]
= 
\left[
\begin{tabular}{c}
1\\
0
\end{tabular}
\right]
= {\top}
\end{align*}

All other logical operators are binary, and hence modelled as rank 3 tensors. To make talking about rank 3 tensors used to model binary operations easier, I will use the following block matrix notation for $2 \times 2 \times 2$ rank 3 tensors $\mathbf{T}$:
\[
\mathbf{T} =
\left[
\begin{tabular}{cc|cc}
$a_1$ & $b_1$ & $a_2$  & $b_2$\\
$c_1$ & $d_1$ & $c_2$  & $d_2$ 
\end{tabular}
\right]
\]
which allows us to express tensor contractions as follows:
\begin{align*}
\mathbf{T} \times \mathbf{v} & =
\left[
\begin{tabular}{cc|cc}
$a_1$ & $b_1$ & $a_2$  & $b_2$\\
$c_1$ & $d_1$ & $c_2$  & $d_2$ 
\end{tabular}
\right]
\left[
\begin{tabular}{c}
$\alpha$\\
$\beta$
\end{tabular}
\right]\\
& =
\left[
\begin{tabular}{cc}
$\alpha \cdot a_1 + \beta \cdot a_2$ & $\alpha \cdot b_1 + \beta \cdot b_2$\\
$\alpha \cdot c_1 + \beta \cdot c_2$ & $\alpha \cdot d_1 + \beta \cdot d_2$
\end{tabular}
\right]
\end{align*}
or more concretely:
{\small
\begin{align*}
& \mathbf{T} \times {\top} =
\left[
\begin{tabular}{cc|cc}
$a_1$ & $b_1$ & $a_2$  & $b_2$\\
$c_1$ & $d_1$ & $c_2$  & $d_2$ 
\end{tabular}
\right]
\left[
\begin{tabular}{c}
$1$\\
$0$
\end{tabular}
\right]
=
\left[
\begin{tabular}{cc}
$a_1$ & $b_1$\\
$c_1$ & $d_1$
\end{tabular}
\right]\\
& \mathbf{T} \times {\bot} =
\left[
\begin{tabular}{cc|cc}
$a_1$ & $b_1$ & $a_2$  & $b_2$\\
$c_1$ & $d_1$ & $c_2$  & $d_2$ 
\end{tabular}
\right]
\left[
\begin{tabular}{c}
$0$\\
$1$
\end{tabular}
\right]
=
\left[
\begin{tabular}{cc}
$a_2$ & $b_2$\\
$c_2$ & $d_2$
\end{tabular}
\right]
\end{align*}}

Using this notation, we can define tensors for the following operations:
\begin{align*}
(\lor) \mapsto \mathbf{T}^{\lor} & =
\left[
\begin{tabular}{cc|cc}
1 & 1 & 1 & 0\\
0 & 0 & 0 & 1
\end{tabular}
\right]\\
 (\land) \mapsto \mathbf{T}^{\land} & =
\left[
\begin{tabular}{cc|cc}
1 & 0 & 0 & 0\\
0 & 1 & 1 & 1
\end{tabular}
\right]\\
 (\to) \mapsto \mathbf{T}^{\to} & =
\left[
\begin{tabular}{cc|cc}
1 & 0 & 1 & 1\\
0 & 1 & 0 & 0
\end{tabular}
\right]
\end{align*}
I leave the trivial proof by exhaustion that these fit the bill to the reader.

It is worth noting here that these tensors preserve normalised probabilities of truth. Let us consider a model such at that described in~\newcite{Coecke2010} which, in lieu of boolean truth values, represents truth value vectors of the form $[\alpha\ \beta]^\top$ where $\alpha + \beta = 1$. Applying the above logical operations to such vectors produces vectors with the same normalisation property. This is due to the fact that the columns of the component matrices are all normalised (i.e.~each column sums to $1$). To give an example with conjunction, let $\mathbf{v} = [\alpha_1\ \beta_1]^\top$ and $\mathbf{w} = [\alpha_2\ \beta_2]^\top$ with $\alpha_1 + \beta_1 = \alpha_2 + \beta_2 = 1$. The conjunction of these vectors is calculated as follows:
\begin{align*}
  &(\mathbf{T}^\land \times \mathbf{v}) \times \mathbf{w} \\
  & \quad = 
  \left[
\begin{tabular}{cc|cc}
1 & 0 & 0 & 0\\
0 & 1 & 1 & 1
\end{tabular}
\right]
\left[
\begin{tabular}{c}
$\alpha_1$\\
$\beta_1$
\end{tabular}
\right]
\left[
\begin{tabular}{c}
$\alpha_2$\\
$\beta_2$
\end{tabular}
\right]\\
& \quad =
\left[
\begin{tabular}{cc}
$\alpha_1$ & $0$\\
$\beta_1$ & $\alpha_1 + \beta_1$ 
\end{tabular}
\right]
\left[
\begin{tabular}{c}
$\alpha_2$\\
$\beta_2$
\end{tabular}
\right]\\
& \quad =
\left[
\begin{tabular}{c}
$\alpha_1 \alpha_2$\\
$\beta_1 \alpha_2 + (\alpha_1 + \beta_1) \beta_2 $ 
\end{tabular}
\right]
\end{align*}
To check that the probabilities are normalised we calculate:
\begin{align*}
& \alpha_1 \alpha_2 + \beta_1 \alpha_2 + (\alpha_1 + \beta_1) \beta_2\\
& = (\alpha_1 + \beta_1) \alpha_2 + (\alpha_1 + \beta_1) \beta_2\\
& = (\alpha_1 + \beta_1)(\alpha_2 + \beta_2) = 1
\end{align*}
We can observe that the resulting probability distribution for truth is still normalised. The same property can be verified for the other connectives, which I leave as an exercise for the reader.


\section{Quantifiers and non-linearity} 
\label{sec:quantifiers_and_non_linearity}

The predicate calculus described up until this point has repeatedly been qualified as `quantifier-free', for the simple reason that quantification cannot be modelled if each application of a predicate or relation immediately yields a truth value. In performing such reductions, we throw away the information required for quantification, namely the information which indicates \emph{which} elements of a domain the predicate holds true or false for. In this section, I present a variant of the predicate calculus developed earlier in this paper which allows us to model simple quantification (i.e.~excluding embedded quantifiers) alongside a tensor-based approach to predicates. However, I will prove that this approach to quantifier modelling relies on non-linear functions, rendering them non-suitable for compositional distributional models relying solely on multilinear maps for composition (or alternatively, rendering such models unsuitable for the modelling of quantifiers by this method).

We saw, in Section~\ref{sec:tensor_based_model_theoretic_predicate_calculi}, that vectors in the semantic space $D$ standing for the logical domain could model logical atoms as well as \emph{sets of atoms}. With this in mind, instead of modelling a predicate $P$ as a truth-function, let us now view it as standing for some function $f_P: \mathcal{P}(\mathcal{D}) \to \mathcal{P}(\mathcal{D})$, defined as:
    \[
    f_P(X) = X \cap M_P
    \]
    where $X$ is a set of domain objects, and $M_P$ is the set modelling the predicate. The tensor form of such a function will be some $\mathbf{T}^{f_P}$in $D \otimes D$. Let this square matrix be a diagonal matrix such that basis weights $c^{T_{f_p}}_{ii} = 1$ if the atom $x$ corresponding to $\mathbf{d}_i$ is in $M_P$ and $0$ otherwise. Through tensor contraction, this tensor maps subsets of $\mathcal{D}$ (elements of $D$) to subsets of $\mathcal{D}$ containing only those objects of the original subset for which $P$ holds (i.e.~yielding another vector in $D$).

To give an example: let us consider a domain with two dogs ($a$ and $b$) and a cat ($c$). One of the dogs ($b$) is brown, as is the cat. Let $S$ be the set of dogs, and $P$ the predicate ``brown''. I represent these statements in the model as follows:
\[
   \mathcal{D} = \{a,\, b,\, c\}\quad S = \{a,\,b\}\quad M_P = \{b,\, c\}
\]
The set of dogs is represented as a vector $\mathbf{S} = [1\ 1\ 0]^\top$ and the predicate `brown' as a tensor in $D \otimes D$:
\[
\mathbf{T}^{P} =
\left[
\begin{tabular}{ccc}
  0 & 0 & 0\\
  0 & 1 & 0\\
  0 & 0 & 1
\end{tabular}
\right]
\]
The set of brown dogs is obtained by computing $f_B(S)$, which distributionally corresponds to applying the tensor $\mathbf{T}^{P}$ to the vector representation of $S$ via tensor contraction, as follows:
\[
\mathbf{T}^{P} \times \mathbf{S} =
\left[
\begin{tabular}{ccc}
  0 & 0 & 0\\
  0 & 1 & 0\\
  0 & 0 & 1
\end{tabular}
\right]
\left[
\begin{tabular}{c}
1\\
1\\
0
\end{tabular}
\right]
=
\left[
\begin{tabular}{c}
0\\
1\\
0
\end{tabular}
\right]
=
\mathbf{b}
\]
The result of this computation shows that the set of brown dogs is the singleton set containing the only brown dog, $b$. As for how logical connectives fit into this picture, in both approaches discussed below, conjunction and disjunction are modelled using set-theoretic intersection and union, which are simply the component-wise $min$ and $max$ functions over vectors, respectively.

Using this new way of modelling predicates as tensors, I turn to the problem of modelling quantification. We begin by putting all predicates in vector form by replacing each instance of the bound variable with a vector $\mathbf{1}$ filled with ones, which extracts the diagonal from the predicate matrix. 

An intuitive way of modelling universal quantification is as follows: expressions of the form ``All $X$s are $Y$s'' are true if and only if $M_X = M_X \cap M_Y$, where $M_X$ and $M_Y$ are the set of $X$s and the set of $Y$s, respectively. Using this, we can define the map $\mathit{forall}$ for distributional universal quantification modelling expressions of the form ``All $X$s are $Y$s'' as follows:
\[
\mathit{forall}(\mathbf{X},\mathbf{Y}) =
\left\{
\begin{tabular}{ll}
$\mathbf{\top}$ & $\quad$ if $\mathbf{X} = min(\mathbf{X},\mathbf{Y})$\\
$\mathbf{\bot}$ & $\quad$ otherwise
\end{tabular}
\right.
\]
To give a short example, the sentence `All Greeks are human' is verified by computing $\mathbf{X} = (\mathbf{M}^{\text{greek}} \times \mathbf{1})$, $\mathbf{Y} = (\mathbf{M}^{\text{human}} \times \mathbf{1})$, and verifying the equality $\mathbf{X} = min(\mathbf{X},\mathbf{Y})$.

Existential statements of the form ``There exists X'' can be modelled using the function $\mathit{exists}$, which tests whether or not $M_X$ is empty, and is defined as follows:
\[
\mathit{exists}(\mathbf{X}) =
\left\{
\begin{tabular}{ll}
$\mathbf{\top}$ & $\quad$ if $|\mathbf{X}| > 0$\\
$\mathbf{\bot}$ & $\quad$ otherwise
\end{tabular}
\right.
\]
To give a short example, the sentence `there exists a brown dog' is verified by computing $\mathbf{X} = (\mathbf{M}^{\text{brown}} \times \mathbf{1}) \cap (\mathbf{M}^{\text{dog}} \times \mathbf{1})$ and verifying whether or not $\mathbf{X}$ is of strictly positive length.

An important point to note here is that neither of these quantification functions are multi-linear maps, since a multilinear map must be linear in all arguments. A counter example for $\mathit{forall}$ is to consider the case where $M_X$ and $M_Y$ are empty, and multiply their vector representations by non-zero scalar weights $\alpha$ and $\beta$.
\begin{align*}
  & \alpha\mathbf{X} = \mathbf{X}\\
  & \beta\mathbf{Y} = \mathbf{Y}\\
  & \mathit{forall}(\alpha\mathbf{X},\beta\mathbf{Y}) = \mathit{forall}(\mathbf{X},\mathbf{Y}) = \mathbf{\top}\\
  & \mathit{forall}(\alpha\mathbf{X}, \beta\mathbf{Y}) \neq \alpha\beta\mathbf{\top}
\end{align*}
I observe that the equations above demonstrate that $\mathit{forall}$ is not a multilinear map.

The proof that $\mathit{exists}$ is not a multilinear map is equally trivial. Assume $M_X$ is an empty set and $\alpha$ is a non-zero scalar weight:
\begin{align*}
  & \alpha \mathbf{X} = \mathbf{X}\\
  & \mathit{exists}(\alpha \mathbf{X}) = \mathit{exists}(\mathbf{X}) = \mathbf{\bot}\\
  & \mathit{exists}(\alpha \mathbf{X}) \neq \alpha \mathbf{\bot}
\end{align*}
It follows that $\mathit{exists}$ is not a multi-linear function.


\section{Conclusions and future work} 
\label{sec:conclusions_and_future_work}

In this paper, I set out to demonstrate that it was possible to replicate most aspects of predicate logic using tensor-based models. I showed that tensors can be constructed from logical models to represent predicates and relations, with vectors encoding elements or sets of elements from the logical domain. I discussed how tensor contraction allows for evaluation of logical expressions encoded as tensors, and that logical connectives can be defined as tensors to form a full quantifier-free predicate calculus. I exposed some of the limitations of this approach when dealing with variables under the scope of quantifiers, and proposed a variant for the tensor representation of predicates which allows us to deal with quantification. Further work on tensor-based modelling of quantifiers should ideally seek to reconcile this work with that of \newcite{Barwise1981}. In this section, I discuss how both of these approaches to predicate modelling can be put into relation, and suggest further work that might be done on this topic, and on the topic of integrating this work into compositional distributional models of semantics.

The first approach to predicate modelling treats predicates as truth functions represented as tensors, while the second treats them as functions from subsets of the domain to subsets of the domain. Yet both representations of predicates contain the same information. Let $\mathbf{M}^P$ and $\mathbf{M'}^P$ be the tensor representations of a predicate $P$ under the first and second approach, respectively. The relation between these representations lies in the equality $diag(\mathbf{p}\mathbf{M}^P) = \mathbf{M'^P}$, where $\mathbf{p}$ is the covector $[1\ 0]$ (and hence $\mathbf{p}\mathbf{M}^P$ yields the first row of $\mathbf{M}^P$). The second row of $\mathbf{M}^P$ being defined in terms of the first, one can also recover $\mathbf{M}^P$ from the diagonal of $\mathbf{M'}^P$.

Furthermore, both approaches deal with separate aspects of predicate logic, namely applying predicates to logical atoms, and applying them to bound variables. With this in mind, it is possible to see how both approaches can be used sequentially by noting that tensor contraction allows for partial application of relations to logical atoms. For example, applying a binary relation to its first argument under the first tensor-based model yields a predicate. Translating this predicate into the second model's form using the equality defined above then permits us to use it in quantified expressions. Using this, we can evaluate expressions of the form ``There exists someone who John loves''. Future work in this area should therefore focus on developing a version of this tensor calculus which permits seamless transition between both tensor formulations of logical predicates.

Finally, this paper aims to provide a starting point for the integration of logical aspects into compositional distributional semantic models. The work presented here serves to illustrate how tensors can simulate logical elements and operations, but does not address (or seek to address) the fact that the vectors and matrices in most compositional distributional semantic models do not cleanly represent elements of a logical domain. However, such distributional representations can arguably be seen as representing the properties objects of a logical domain hold in a corpus: for example the similar distributions of `car' and `automobile' could serve to indicate that these concepts are co-extensive. This suggests two directions research based on this paper could take. One could use the hypothesis that similar vectors indicate co-extensive concepts to infer a (probabilistic) logical domain and set of predicates, and use the methods described above without modification; alternatively one could use the form of the logical operations and predicate tensors described in this paper as a basis for a higher-dimensional predicate calculus, and investigate how such higher-dimensional `logical' operations and elements could be defined or learned. Either way, the problem of reconciling the fuzzy `messiness' of distributional models with the sharp `cleanliness' of logic is a difficult problem, but I hope to have demonstrated in this paper that a small step has been made in the right direction.


\section*{Acknowledgments}

Thanks to Ond\v{r}ej Ryp\'{a}\v{c}ek, Nal Kalchbrenner and Karl Moritz Hermann for their helpful comments during discussions surrounding this paper. This work is supported by EPSRC Project \texttt{EP/I03808X/1}.


\end{document}